\documentclass[11pt,reqno]{amsart}
\usepackage{amsmath,amssymb,amsthm}
\usepackage{mathrsfs}
\usepackage{enumitem}
\usepackage[margin=1in]{geometry}
\usepackage[hidelinks]{hyperref}

\DeclareMathOperator{\supp}{supp}
\DeclareMathOperator{\Lie}{Lie}
\newcommand{\R}{\mathbb{R}}

\theoremstyle{plain}
\newtheorem{theorem}{Theorem}[section]
\newtheorem{proposition}[theorem]{Proposition}
\newtheorem{corollary}[theorem]{Corollary}

\theoremstyle{definition}
\newtheorem{definition}[theorem]{Definition}
\newtheorem{example}[theorem]{Example}
\theoremstyle{remark}
\newtheorem{remark}[theorem]{Remark}
\newcommand{\step}[1]{\medskip\par\noindent\textit{#1}\quad}

\title[Liouville theorem for L\'evy operators on nilpotent and compact Lie groups]
{A Liouville theorem for L\'evy operators\\ on nilpotent and compact Lie groups}

\author{Samuel Killick}
\address{School of Mathematical and Physical Sciences, University of Technology Sydney (UTS)}
\email{samuel.killick@uts.edu.au}

\author{David Lee}
\address{School of Mathematical and Physical Sciences, University of Technology Sydney (UTS)}
\email{davidchanwoo.lee@uts.edu.au}

\subjclass[2020]{35A08, 35R11, 42B37, 60B15, 60G51}
\keywords{Liouville theorem, nonlocal operator, L\'evy process, nilpotent Lie group, Choquet--Deny}

\begin{document}

\begin{abstract}
Let $\mathcal{L}$ be a left-invariant L\'evy operator, with characteristics $(a,b,\mu)$, on a connected
    and simply connected Lie group $G$ that is either compact or nilpotent with $[G,[G,G]]$ relatively
    compact (equivalently, of step at most two). We prove that every bounded distributional solution of $\mathcal{L}u=0$ is right-invariant under
    the closed subgroup $P=\overline{\langle\exp\mathfrak h_0\cup\supp\mu\rangle}$, where
    $\mathfrak h_0=\Lie\langle b,\operatorname{ran}(a)\rangle$, and is constant if and only if $P=G$. In the
    compact case we give three proofs, one through the It\^o--Kawada equidistribution theorem, one by
    propagation of the maximum and one by an application of the Peter-Weyl theorem, while the nilpotent case rests on the Choquet--Deny theorem of B\'enard and Frisch.
\end{abstract}

\maketitle

\section{Introduction}\label{sec:intro}

The classical theorem of Liouville states that a bounded harmonic function on
$\R^{d}$, that is a bounded function $u$ with $\Delta u=0$, must be constant. Over the last
decade this principle has been revisited in the setting of \emph{nonlocal} operators.
In the case of the fractional Laplacian $(-\Delta)^{s}$, $s\in(0,1)$, the Liouville property
holds classically as a consequence of the Harnack inequality. Fall \cite{Fall2016} proved that every entire $s$-harmonic
function is \emph{affine}, and is constant precisely when $s\in(0,\tfrac12]$. Moreover, he
showed that if $(-\Delta)^{s}u$ equals a polynomial then that polynomial must
vanish and $u$ is again affine.

For the more general setting of L\'evy operators, which can be expressed as
\begin{equation}\label{eq:euclidean-levy}
    \mathcal{L} u(x)=\operatorname{tr}\!\big(a\,D^{2}u(x)\big)+b\cdot\nabla u(x)
    +\int_{\R^{d}\setminus\{0\}}\!\Big(u(x+z)-u(x)-z\cdot\nabla u(x)\mathbf 1_{|z|\le1}\Big)\mu(dz),
\end{equation}
with L\'evy triplet $(a,b,\mu)$, and the question becomes much more subtle. Concretely,
one can ask for which triplets is every bounded solution of $\mathcal{L} u=0$ constant?

This question was first answered by Alibaud, del~Teso, Endal and Jakobsen \cite{ADEJ}. More specifically,
they characterized the Liouville property for L\'evy operators by using a blend of the maximum principle and propagation of a maximum
across the support of the jump measure.
This allowed them to obtain a sharp condition on
the triplet $(a,b,\mu)$. They showed that bounded solutions are periodic with respect to an
explicit closed subgroup of $\R^{d}$ that can be characterized by the triplet $(a,b,\mu)$. In particular, when the group is all
of $\R^{d}$ they show that solutions are necessarily constant.

Shortly afterwards, Knopova and Schilling \cite{KnopovaSchilling}, K\"uhn \cite{Kuhn} and then Berger and Schilling \cite{BergerSchilling} characterized the same property (through different assumptions on the solution) through Fourier analysis and distribution theory, recasting it in terms of the characteristic exponent (symbol) $\psi$ of the operator. Via a division-type argument, they show that the Fourier transform of a $\mathcal{L}$-harmonic function is a distribution supported on the zero set $\{\psi=0\}$, so when that set reduces to $\{0\}$ the solution must be constant, whereas any further zero $\xi_{0}\neq0$ yields the bounded nonconstant solution $x\mapsto\cos(\xi_{0}\cdot x)$. Berger, Schilling and Shargorodsky \cite{BSS} later re-proved the criterion for a general class of Fourier multipliers by a purely analytic argument, linking it to the coupling property of the associated semigroup. More recently, \cite{BSSS} extends the characterization from bounded to sub-exponentially growing solutions. A related division-type classification, valid without any boundedness or decay hypothesis on the solution, is obtained in \cite{LeeDivision} for symbols vanishing to arbitrary finite order on their zero set, and the same method yields a Liouville theorem for the logarithmic Laplacian \cite{LeeLogLaplacian}.

Grzywny and Kwa\'snicki \cite{GrzywnyKwasnicki} go further, handling positive and
polynomially bounded $\mathcal{L}$-harmonic functions under the same non-degeneracy
$\{\psi=0\}=\{0\}$. Their main tool is a representation theorem of Deny type. We call an exponential
$e_\lambda(x)=e^{\lambda\cdot x}$ \emph{$\mathcal{L}$-harmonic} when $\mathcal{L} e_\lambda=0$, and
write $\Lambda$ for the set of $\lambda$ for which this holds. They show that every
nonnegative solution is built out of these, in that any $h\ge0$ with $\mathcal{L} h=0$ can be written as
$h(x)=\int_{\Lambda}e^{\lambda\cdot x}\,\sigma(d\lambda)$ for a unique nonnegative
measure $\sigma$ on $\Lambda$. Since the only bounded
exponentials are constant, the Liouville property drops out as a special case. For
polynomial growth they show, assuming $\psi$ is smooth away from the origin (or suitable
moment bounds on the L\'evy measure), that $\mathcal{L}$-harmonic functions of polynomial growth
are polynomials with a particular degree.

Nevertheless, all of the work mentioned above takes place on $\R^{d}$, where the Euclidean
Fourier transform diagonalises the operator and translation commutes with
differentiation.

\textbf{This article aims to prove a Liouville theorem in the case where one works on a compact or nilpotent Lie group (under appropriate assumptions). We also compare 3 methods which have been used in the Euclidean case in this setting. }

We outline each approach in reference to the existing approaches within the Euclidean setting.

\begin{itemize}

    \item \emph{A propagation-of-maximum proof.} Alibaud, del~Teso, Endal and Jakobsen
          \cite{ADEJ} proved the Euclidean Liouville property via the positive
          maximum principle, by moving a maximum across the jump support and along the diffusion and
          drift. The same argument runs on a compact group.

    \item \emph{A Choquet--Deny proof.} Both Grzywny and Kwa\'snicki \cite{GrzywnyKwasnicki} and Berger and Schilling \cite{BergerSchilling} reduce the Liouville property for a L\'evy
          generator to the convolution equation involving the transition probability measure associated with the L\'evy operator. It is from here that Berger and Schilling are able to prove their Liouville theorem via Fourier analytic methods. We apply a similar proof, in the compact setting, but we deviate from their approach by using the \emph{resolvent measure} and the \emph{It\^o--Kawada theorem}.

    \item \emph{A harmonic analysis proof.} Over $\R^{d}$, Knopova--Schilling and
          Berger--Schilling \cite{KnopovaSchilling,BergerSchilling} derive the Liouville property off
          the characteristic exponent, Berger, Schilling and Shargorodsky \cite{BSS} handle general
          Fourier multipliers. On a compact group, Peter--Weyl takes the place of the Fourier
          transform. Expanding $u$ over the unitary dual, $\mathcal{L}$ acts on each isotypic component as an
          operator-valued symbol $\psi(\pi)$, and $\psi(\pi)$ is invertible for every non-trivial
          $\pi$, so only the trivial component survives and $u$ is constant. The scalar condition
          $\{\psi=0\}=\{0\}$ becomes the invertibility of $\psi(\pi)$ for all non-trivial $\pi$.
\end{itemize}

The \emph{nilpotent} case requires more consideration. In this setting we can still establish a Liouville theorem (under additional conditions on the group) via a \emph{Choquet--Deny} argument.

The Choquet--Deny property has a long history. A probability measure $\nu$ on a locally compact group $H$ is \emph{adapted} if its support generates a dense subgroup, and $H$ is a \emph{Choquet--Deny group} if for every adapted $\nu$ the only bounded $\nu$-harmonic functions are constant. Choquet and Deny~\cite{ChoquetDeny} established this property for all abelian locally compact groups, and Dynkin and Malyutov~\cite{DynkinMalyutov} proved it for finitely generated nilpotent groups. For connected Lie groups, Azencott~\cite{Azencott} and Guivarc'h~\cite{Guivarch} identified those that are Choquet--Deny for spread-out measures as the groups of polynomial growth, and Guivarc'h recorded the question, due to Furstenberg, of whether every nilpotent Lie group is Choquet--Deny. In the discrete setting the countable Choquet--Deny groups were characterized by Frisch, Hartman, Tamuz and Vahidi Ferdowsi~\cite{FHTV} through the absence of quotients with the infinite conjugacy class property. For general measures on nilpotent Lie groups the question proved delicate, an earlier claim of Raugi~\cite{Raugi} having been withdrawn, and it was settled only recently by B\'enard and Frisch~\cite{BenardFrisch}, who proved that a connected nilpotent Lie group is Choquet--Deny if and only if its iterated commutator subgroup $[G,[G,G]]$ is relatively compact. It is this result, stated as Theorem~\ref{thm:benard-frisch}, that we use to treat the nilpotent case.

Before we present our main results, we introduce some notations and conventions.

\subsection{Notation and Conventions}\label{subsec:prelim}

From this point onwards, we will use the following notations and conventions. We denote $G$ to be a connected, simply connected, unimodular Lie group (of dimension $d$) with identity $e$, Lie algebra $\Lie(G)$ and a fixed basis of left-invariant vector fields $X_1,\dots,X_{d}$. Moreover, we also denote $\exp:\Lie(G)\rightarrow G$ to be the exponential map. With this choice, we fix smooth canonical coordinates $x_{1},\dots,x_{d}$ such that
\begin{itemize}
    \item $x_{i}(e)=0$, for $i=1,\dots,d $,
    \item $X_{i}x_{j}(e)=\delta_{ij}$, $i,j=1,\dots,d$ where $\delta_{ij}$ is the Kronecker delta,
    \item $x_{i}(\tau^{-1})=-x_{i}(\tau)$ for all $\tau \in G$.
\end{itemize}
We denote by $d\tau$ the Haar measure on $G$. Due to the unimodular assumption, the Haar measure has the following properties.
\begin{itemize}
    \item Bi-invariance $\displaystyle \int_{G}f(\tau)\,d\tau=\int_{G}f(g\star\tau)\,d\tau=\int_{G}f(\tau \star g)\,d\tau,$ for $g\in G$.
    \item Inversion-invariance $\displaystyle \int_{G}f(\tau)\,d\tau=\int_{G}f(\tau^{-1})\,d\tau$.
\end{itemize}
We fix a relatively compact open neighborhood  $B_0$ containing $e$ on which the \emph{exponential map} $\exp:\Lie(G)\rightarrow G$ is a diffeomorphism onto its image and let $y=\log:B_0\rightarrow \Lie(G)$ be the \emph{canonical coordinates of the first kind}. From this, we can fix a family of smooth bump functions of compact support $x_1,\dots,x_d$ such that
\begin{equation*}
    x_i=y_i, \quad \text{on a neighborhood of $e$,}\quad \supp(x_i)\subset B_0.
\end{equation*}
We denote $L_{h}f=f(h\star\cdot)$ and $R_{h}f=f(\cdot\star h)$ to be the left and right
translations (respectively). We use two convolutions. For functions $\psi,f$ on $G$ we write
the \emph{left convolution}
\begin{equation*}
    (\psi*f)(g):=\int_{G}\psi(\tau)\,f(\tau^{-1}\star g)\,d\tau,\qquad g\in G,
\end{equation*}
whenever the integral is well-defined. This is the convolution we use to mollify functions.
For a bounded measurable function $f$ and a finite Borel measure $\nu$ on $G$ we write the
\emph{convolution of $f$ with $\nu$}
\begin{equation*}
    (f*\nu)(g):=\int_{G}f(g\star h)\,\nu(dh),\qquad g\in G,
\end{equation*}
which is the convolution appearing in the transition semigroup. Finally, for finite Borel
measures $\mu,\nu$ on $G$ their convolution $\mu*\nu$ is the push-forward of $\mu\otimes\nu$
under the group product,
\begin{equation*}
    \int_{G}f\,d(\mu*\nu)=\int_{G}\int_{G}f(h\star h')\,\mu(dh)\,\nu(dh'),
\end{equation*}
so that $(f*\nu)*\mu=f*(\mu*\nu)$ for every bounded measurable $f$.
We also adopt the following conventions for function spaces.

\begin{itemize}
    \item $C(G)$ continuous functions on $G$,
    \item $C_b(G)$ continuous and bounded functions on $G$,
    \item $C_{c}^\infty(G)=\mathcal{D}(G)$ smooth functions with compact support on $G$ (test functions),
    \item $\mathcal{D}'(G)$ space of distributions on $G$, the dual of $\mathcal{D}(G)$.
\end{itemize}

We now introduce the operators of interest in this note.

\begin{definition}\label{def:levy_operator}
    A left-invariant L\'evy operator on $G$ is a linear map $\mathcal{L}:C^{\infty}_{c}(G)\rightarrow C(G)$ of the form
    \begin{equation}\label{eq:generator}
        \mathcal{L}f(g)=\sum_{i,j=1}^{d}a_{ij}X_{i}X_{j}f(g)+\sum_{i=1}^{d}b_{i}X_{i}f(g)+\int_{G}\left (f(g\star\tau)-f(g)-\sum_{i=1}^{d}x_{i}(\tau)X_{i}f(g) \right )\,\mu(d\tau),
    \end{equation}
    for $g\in G$, where
    \begin{itemize}
        \item $a=(a_{ij})$ is a symmetric, positive semidefinite matrix.
        \item $b=(b_i)\in \R^d$.
        \item $\mu$ is a L\'evy measure. More specifically, $\mu$ is a nonnegative Borel measure on $G\setminus\{e\}$ such that
              \begin{equation*}
                  \int_{V} \sum_{i=1}^{d}x_i(\tau)^2 \mu(d\tau)<\infty, \quad \text{and }\quad \mu(V^{c})<\infty,
              \end{equation*}
              for every neighborhood $V$ which contains $e$.
    \end{itemize}
    We call $(a,b,\mu)$ the characteristics of $\mathcal{L}$.
\end{definition}

We remark that left-invariant L\'evy operators are precisely the left-invariant operators (with a zero killing term) which satisfy the positive maximum principle via the Courr\'ege theorem due to Applebaum and Ngan~\cite{MR4165724,MR4072488}. Moreover, via Hunt's theorem on Lie groups~\cite{MR79232} (see also~\cite[\S 3.1]{MR2060091}), we know that there exists a unique convolution semigroup $(\rho_{t})_{t\geq 0}$ of probability measures on $G$ which satisfies
\begin{itemize}
    \item $\rho_0=\delta_e$,
    \item $\rho_{s+t}=\rho_{s}*\rho_{t}$, for all $s,t \geq 0$,
    \item $\rho_{t}\rightarrow \delta_e$ weakly as $t\rightarrow 0$.
\end{itemize}

The corresponding Feller semigroup (on $C_0(G)$) is given by $T_{t}f=f*\rho_{t}$ and its generator (restricted to the core $C^\infty_{c}(G)$) is given by
\begin{equation*}
    \frac{d}{dt}T_{t}f\big |_{t=0} = \mathcal{L}f, \quad \text{for $f\in C^\infty_{c}(G)$. }
\end{equation*}

We will write $(X_{t})_{t\geq 0}$ to be the corresponding L\'evy process in $G$ such that
\begin{itemize}
    \item $X_{0}=g\in G$,
    \item $\mathbb{P}(X_{t}\in A)=\rho_{t}(A)$, for Borel sets $A$,
    \item $\mathbb{E}(f(X_{t})|X_0=g)=(f*\rho_{t})(g)$ for $g\in G$ and $f\in C_b(G)$.
\end{itemize}

We use the \emph{resolvent measure} (which is a probability measure) given by
\begin{equation*}
    \rho(A):=\int_{0}^\infty e^{-t}\rho_{t}(A) \,dt, \quad \text{for Borel measurable sets $A\subset G$. }
\end{equation*}

\subsection{Main result}\label{sec:main}

For the main result, we wish to characterize bounded solutions of the following problem
\begin{equation}\label{eqn:main_eqn}
    \mathcal{L}u =0, \quad \text{on $G$,}
\end{equation}
where $\mathcal{L}$ is a left-invariant L\'evy operator with triplet $(a,b,\mu)$. Before we state our result we introduce the appropriate notions of solution corresponding to~\eqref{eqn:main_eqn}.

\begin{definition}\label{def:notion_of_solution}
    Let $u\in L^\infty(G)$ and $\mathcal{L}$ be a left-invariant L\'evy operator. We define the distribution $\mathcal{L}u\in \mathcal{D}'(G)$ by
    \begin{equation*}
        \langle \mathcal{L}u,\varphi\rangle_{\mathcal{D}'(G),\mathcal{D}(G)}:=\int_{G}u\left (\mathcal{L}^\dagger\varphi \right )\,d\tau, \quad \text{for all $\varphi \in \mathcal{D}(G)$,}
    \end{equation*}
    where $\mathcal{L}^\dagger$ is the Haar-adjoint of $\mathcal{L}$.
    We say that $u\in L^\infty(G)$ is a bounded distributional solution of~\eqref{eqn:main_eqn} if we have that
    \begin{equation*}
        \langle \mathcal{L}u,\varphi \rangle_{\mathcal{D}'(G),\mathcal{D}(G)}=0, \quad \text{for all $\varphi \in \mathcal{D}(G)$.}
    \end{equation*}
\end{definition}

\begin{remark}
    We show that the notion of bounded distributional solution is indeed valid. It can be verified by direct computation that the Haar adjoint of $\mathcal{L}$ is again a left-invariant L\'evy operator. Namely, if the L\'evy triplet corresponding to $\mathcal{L}$ is given by $(a,b,\mu)$ then the L\'evy triplet corresponding to $\mathcal{L}^\dagger$  is given by $(a,-b, \hat{\mu})$ where $\hat{\mu}(A)=\mu(A^{-1})$. Here $A^{-1}$ consists of the group inverses of $A$. Note that this critically requires the unimodular assumption of $G$. This gives us that the adjoint is given by
    \begin{equation*}
        \mathcal{L}^{\dagger}f(g)=\sum_{i,j=1}^{d}a_{ij}X_{i}X_{j}f(g)-\sum_{i=1}^{d}b_i X_{i}f(g)+\int_{G}\left ( f(g \star\tau^{-1})-f(g)+\sum_{i=1}^{d}x_{i}(\tau)X_{i}f(g)\right )\mu(d\tau),
    \end{equation*}
    for $g\in G$. Hence to show that the notion of solution is appropriate it is sufficient to demonstrate that $\mathcal{L}\psi\in L^1(G)$ for any L\'evy operator $\mathcal{L}$ and $\psi\in \mathcal{D}(G)$. From an application of the second order Taylor expansion

    \begin{equation*}
        \left | \psi(g\star h)-\psi(g )-\sum_{i}x_{i}(h)X_{i}\psi(g)\right |\leq \frac{1}{2}\max_{i,j}\|X_{i}X_{j}\psi\|_{\infty}\sum_{i}x_{i}(h)^2,\quad \text{for $g, h \in G$, }
    \end{equation*}
    we have that $\mathcal{L}\psi\in L^1(G)$ for any $\psi \in \mathcal{D}(G)$.
\end{remark}
Before we state the main results of this article we introduce the following sets
\begin{equation*}
    P=\overline{\big\langle\,\exp\mathfrak h_0\cup\supp\mu\,\big\rangle},
    \qquad \mathfrak h_0=\Lie\langle b,\operatorname{ran}(a)\rangle .
\end{equation*}

\begin{theorem}\label{thm:main}
    Let $G$ be a connected, simply connected Lie group which is either
    \begin{itemize}
        \item compact,
        \item nilpotent with the property that $[G,[G,G]]$ is relatively compact; as $G$ is simply connected, this is equivalent to $G$ being nilpotent of step at most two, that is, abelian or $2$-step (see Remark~\ref{rem:step-two}).
    \end{itemize}
    Let $u$ be a bounded distributional solution of~\eqref{eqn:main_eqn}, see Definition~\ref{def:notion_of_solution}, for a L\'evy operator with triplet $(a,b,\mu)$. Then we have that $u$ is right-invariant, namely
    \begin{equation*}
        u(g\star p)=u(g), \quad \text{for almost every $g\in G$ and $p\in P$.}
    \end{equation*}
    Moreover, $u$ is constant if $P=G$, while if $P\neq G$ there exist nonconstant bounded distributional solutions.
\end{theorem}

The above result is essentially a Lie group analogue of the result of Alibaud, del~Teso, Endal and Jakobsen~\cite{MR4149690}, however, the proof we use does not rely on their \emph{propagation-of-maximum} approach. Instead, we demonstrate that it is an immediate corollary of a Choquet--Deny type theorem, due to a recent result of B\'enard and Frisch~\cite{BenardFrisch}, see Theorem~\ref{thm:benard-frisch}. We remark that one can establish when a solution is necessarily constant. However, when $P\ne G$, showing that $u$ is right-invariant with respect to $P$ does not quite work when the Lie group is non-commutative. We present this in more detail in Section~\ref{sec:propagation}.

\begin{remark}\label{rem:step-two}
    The nilpotent hypothesis of Theorem~\ref{thm:main} is more restrictive than it first appears. When $G$
    is simply connected and nilpotent, each term of its lower central series is a closed, simply connected
    subgroup, hence isomorphic to some $\R^k$, so $[G,[G,G]]$ is relatively compact only when it is trivial.
    The simply connected nilpotent groups covered by Theorem~\ref{thm:main} are therefore exactly those of
    step at most two. These include the abelian groups $\R^d$, the Heisenberg groups $\mathbb H^n$, and more
    generally all Heisenberg-type (Kaplan) groups and free two-step nilpotent groups, together with direct
    products such as $\R^k\times\mathbb H^n$. A simply connected nilpotent group of step three or higher,
    such as the group of $4\times4$ upper unitriangular matrices, has non-compact iterated commutator
    $[G,[G,G]]$ and is not covered.
\end{remark}

\subsection{Examples}

We first check that Theorem~\ref{thm:main} recovers the known Euclidean theorem, and then illustrate the
genuinely non-commutative behavior on $SU(2)$.

\begin{example}\label{ex:euclidean}
    Take $G=\R^d$ with its additive group structure. It is abelian, hence nilpotent of step one, connected
    and simply connected, and $[G,[G,G]]=\{0\}$ is trivially relatively compact, so $\R^d$ falls under the
    nilpotent case of Theorem~\ref{thm:main}. The left-invariant fields are the commuting partial
    derivatives $X_i=\partial_{x_i}$, the exponential map is the identity, and a left-invariant L\'evy
    operator is the classical L\'evy operator~\eqref{eq:euclidean-levy}.

    Because the group is abelian, $\mathfrak h_0=\Lie\langle b,\operatorname{ran}(a)\rangle$ is simply the
    linear subspace $V=\operatorname{span}\{b\}+\operatorname{ran}(a)$, and $\exp\mathfrak h_0=V$. Hence
    \begin{equation*}
        P=\overline{\langle\,V\cup\supp\mu\,\rangle}=\overline{V+\langle\supp\mu\rangle},
    \end{equation*}
    the closed subgroup of $\R^d$ generated by the subspace $V$ (spanned by the drift and the range of the
    diffusion) together with the jump vectors. Theorem~\ref{thm:main} then states that every bounded
    solution of $\mathcal{L}u=0$ is $P$-periodic, that is $u(x+p)=u(x)$ for all $p\in P$, and is constant if
    and only if $P=\R^d$. This is precisely the characterization of Alibaud, del~Teso, Endal and
    Jakobsen~\cite{MR4149690,ADEJ}, with $P$ their subgroup of periods.

    The two regimes appear in the following concrete choices.
    \begin{itemize}
        \item \emph{Fractional Laplacian.} For $\mathcal{L}=-(-\Delta)^s$ with $s\in(0,1)$ one has $a=0$,
              $b=0$ and $\mu(dz)=c_{d,s}\,|z|^{-d-2s}\,dz$, whose support is all of $\R^d\setminus\{0\}$.
              Thus $\langle\supp\mu\rangle$ is dense, $P=\R^d$, and every bounded solution is constant. This
              is the classical Liouville theorem for the fractional Laplacian.
        \item \emph{A proper jump lattice.} On $\R$ take the pure-jump operator
              $\mathcal{L}u(x)=u(x+1)+u(x-1)-2u(x)$, that is $a=0$, $b=0$ and $\mu=\delta_1+\delta_{-1}$. Here
              $V=\{0\}$ and $\langle\supp\mu\rangle=\mathbb{Z}$, so $P=\mathbb{Z}\neq\R$. A bounded $u$ solves
              $\mathcal{L}u=0$ exactly when it is $1$-periodic, for instance $u(x)=\cos(2\pi x)$, so
              nonconstant bounded solutions exist, in agreement with $P\neq\R$.
    \end{itemize}
\end{example}

We illustrate Theorem~\ref{thm:main} on the simplest genuinely non-commutative compact group, $G=SU(2)$,
which is compact, connected and simply connected, hence covered by the compact case. Fix a basis
$X_1,X_2,X_3$ of $\Lie(SU(2))$ with
\begin{equation*}
    [X_1,X_2]=X_3,\qquad [X_2,X_3]=X_1,\qquad [X_3,X_1]=X_2,
\end{equation*}
and let $T=\{\exp(\theta X_3):\theta\in\R\}$ be the maximal torus generated by $X_3$, so that
$\exp(\theta X_3)=\operatorname{diag}(e^{i\theta/2},e^{-i\theta/2})$. The two examples below form a matched
pair, built on the same group and the same kind of operator and differing only in how many diffusion
directions are switched on. In the first the characteristics generate only $T$, so $P=T\neq G$ and the
Liouville property fails, while in the second one further direction makes them bracket-generate, so $P=G$
and it holds.

\begin{example}[A compact, non-commutative operator that is \emph{not} Liouville]\label{ex:su2}
    Take the single degenerate diffusion $\mathcal{L}=X_3^2$, that is $\operatorname{ran}(a)=\R X_3$, $b=0$
    and $\mu=0$. The only direction present is $X_3$, and so
    \begin{equation*}
        \mathfrak h_0=\Lie\langle\operatorname{ran}(a)\rangle=\R X_3,\qquad
        P=\overline{\langle\exp\R X_3\rangle}=T\neq G .
    \end{equation*}
    Theorem~\ref{thm:main} says that every bounded solution of $X_3^2u=0$ is right-invariant under $T$, and
    since $P=T\neq G$ it need not be constant. Indeed it is not. Writing $g=(g_{ij})\in SU(2)$, set
    \begin{equation*}
        u(g)=|g_{11}|^2 .
    \end{equation*}
    Right multiplication by $\exp(\theta X_3)=\operatorname{diag}(e^{i\theta/2},e^{-i\theta/2})$ multiplies the
    first column of $g$ by the phase $e^{i\theta/2}$ and so leaves $|g_{11}|^2$ unchanged. Hence $u$ is
    right-$T$-invariant, which gives $X_3u=0$ and so $X_3^2u=0$, while $u$ is plainly nonconstant. In fact a
    bounded function solves $X_3^2u=0$ exactly when it is right-$T$-invariant, that is, when it descends to
    the quotient $SU(2)/T\cong S^2$. The bounded solutions are therefore all bounded functions on the
    $2$-sphere, a whole function space rather than the constants alone. This is the non-commutative analogue
    of a bounded harmonic function that depends on only some of the variables.
\end{example}

\begin{example}[Switching the Liouville property back on]
    Keep $G=SU(2)$ and add one further diffusion direction, $\mathcal{L}=X_1^2+X_3^2$, so
    $\operatorname{ran}(a)=\R X_1\oplus\R X_3$. The bracket $[X_3,X_1]=X_2$ recovers the missing direction, so
    $X_1$ and $X_3$ now bracket-generate the whole Lie algebra and
    \begin{equation*}
        \mathfrak h_0=\Lie\langle X_1,X_3\rangle=\Lie(SU(2)),\qquad
        P=\overline{\langle\exp\Lie(SU(2))\rangle}=SU(2)=G .
    \end{equation*}
    This time $P=G$, and so Theorem~\ref{thm:main} forces every bounded solution to be constant. Placing the
    two side by side, the bounded solutions of $X_3^2u=0$ fill $L^\infty(S^2)$, whereas adding the one
    bracket-generating direction $X_1$ collapses them to the constants. The obstruction in the first example
    is thus exactly the failure of $P$ to fill out $G$.
\end{example}

\subsection{Outline of the article}

The remainder of the article is organized as follows. In Section~\ref{sec:reduction} we reduce the problem to smooth solutions of $\mathcal{L}u=0$, taken in the pointwise sense, by a mollification argument adapted from the Euclidean work of Alibaud, del~Teso, Endal and Jakobsen~\cite{MR4149690}. Section~\ref{sec:choquet} proves Theorem~\ref{thm:main} in both the compact and the nilpotent cases: we recast the equation as the Choquet--Deny convolution equation $u=u*\rho$ for the resolvent measure $\rho$, apply the It\^o--Kawada theorem in the compact case and the theorem of B\'enard and Frisch~\cite{BenardFrisch} in the nilpotent case, and identify the subgroup $P$ in terms of the characteristics $(a,b,\mu)$. In Section~\ref{sec:propagation} we revisit the compact case by the propagation-of-maxima method of~\cite{MR4149690}, which recovers the part of Theorem~\ref{thm:main} asserting that $u$ is constant if and only if $P=G$ (Theorem~\ref{thm:compact-propagation}), and we explain why it does not yield the full right-invariance on a non-commutative or non-compact group. Finally, Section~\ref{sec:harmonic} gives a Fourier-analytic proof in the compact case, based on the Peter--Weyl theorem, and records some observations on a Fourier-analytic approach to the nilpotent case.

\section{Reduction to smooth solutions}\label{sec:reduction}
Before we commence on proving Theorem~\ref{thm:main} directly, we demonstrate that we can reduce to studying smooth solutions which solve~\eqref{eqn:main_eqn} in a pointwise sense. This is similar in spirit to what was considered in the Euclidean case by Alibaud, del~Teso, Endal and Jakobsen~\cite{MR4149690}.

The goal of this section is to prove the following proposition which is the Lie group analogue of~\cite[Lemma 2.4]{MR4149690}. This characterization will be useful for when we consider the approaches via \emph{propagation-via-maximum} and \emph{Choquet--Deny}.

\begin{proposition}\label{prop:smooth_soln}
    Let $u\in L^\infty(G)$ and let $\{\psi_{\epsilon}\}_{\epsilon>0}
        \subset C_c^\infty(G)$ be an approximate identity, and set
    $u_{\epsilon}:=\psi_{\epsilon}*u$. Then the following hold.
    \begin{itemize}
        \item[1.] $\mathcal{L}u=0$ in $\mathcal{D}'(G)$ if and only if $\mathcal{L}u_{\epsilon}=0$
              pointwise for every $\epsilon>0$,
        \item[2.] $u$ is a.e.\ $S$-periodic (i.e.\ $R_{s}u=u$ for almost every $s\in S$) if and only if
              $R_{s}u_{\epsilon}=u_{\epsilon}$ for every $s\in S$ and every
              $\epsilon>0$.
    \end{itemize}
\end{proposition}

The above proposition effectively allows us to consider solutions that are smooth with bounded derivatives. Moreover, under this assumption and It\^o's formula~\cite{MR2060091}, we have the following identity which can be viewed as an alternative to~\eqref{eqn:main_eqn}.
\begin{corollary}\label{cor:choquet-deny}
    Let $(\rho_{t})_{t\geq 0}$ be the corresponding convolution semigroup associated with $\mathcal{L}$. Suppose that $u\in C^\infty_{b}(G)$ and $\mathcal{L}u=0$. Then, we have that
    \begin{equation}\label{eqn:convolution_eqn1}
        u=u*\rho_{t}, \quad\text{for all $t\geq 0$. }
    \end{equation}
\end{corollary}

Before we prove Proposition~\ref{prop:smooth_soln}, we remind the reader that $G$ is unimodular since $G$ is compact or nilpotent. So the modular function is $\Delta\equiv 1$. We
write $\tilde\psi(h):=\psi(h^{-1})$ for the reflection and record three
consequences of unimodularity that we use repeatedly.
\begin{itemize}
    \item[(i)] \emph{Adjoint of left convolution.} For $\psi\in C_c^\infty(G)$,
          $w\in L^\infty(G)$ and $\varphi\in L^1(G)$,
          \begin{equation}\label{eq:conv-adjoint}
              \int_{G}(\psi * w)\,\varphi \, dh
              = \int_{G} w\,(\tilde\psi*\varphi) \, dh ,
          \end{equation}
          which is Fubini together with the substitution $h\mapsto h^{-1}$ under the
          (two-sided) Haar measure, here $\psi^{\vee}=\tilde\psi$ because $\Delta\equiv1$.
    \item[(ii)] \emph{Left-invariance of the adjoint.} Since $G$ is unimodular,
          left-invariant vector fields are formally skew-adjoint, so the formal adjoint
          $\mathcal{L}^{\dagger}$ of the left-invariant generator $\mathcal{L}$ is again left-invariant. In
          particular $\mathcal{L}^{\dagger}$ commutes with left convolution,
          \begin{equation*}
              \mathcal{L}^{\dagger}(\varphi*\psi)=\varphi*(\mathcal{L}^{\dagger}\psi),
              \qquad \varphi,\psi\in\mathcal{D}(G).
          \end{equation*}
    \item[(iii)] \emph{Left convolution commutes with right translation:}
          $R_{s}(\psi*w)=\psi*(R_{s}w)$ for $s\in G$.
\end{itemize}

We can now prove Proposition~\ref{prop:smooth_soln}.
\begin{proof}[Proof of Proposition~\ref{prop:smooth_soln}]
    First observe that each $u_{\epsilon}=\psi_{\epsilon}*u$ is smooth and bounded, so $\mathcal{L}u_{\epsilon}$ and $R_{s}u_{\epsilon}-u_{\epsilon}$ are continuous functions.

    \begin{itemize}
        \item[1.] We consider the forward direction. Assume that
              \begin{equation*}
                  \int_{G}u\mathcal{L}^{\dagger}\varphi\, d\tau=0, \quad \text{for all $\varphi\in
                          \mathcal{D}(G)$. }
              \end{equation*}
              Replacing $\varphi$ by $\tilde\psi_{\epsilon}*\varphi\in \mathcal{D}(G)$, then using
              the left-invariance of $\mathcal{L}^{\dagger}$ and the transpose
              identity for left convolution, we obtain
              \begin{equation*}
                  0=\int_{G}u\,\mathcal{L}^{\dagger}\!\bigl(\tilde\psi_{\epsilon}*\varphi\bigr)\,d\tau
                  =\int_{G}u\,\bigl(\tilde\psi_{\epsilon}*\mathcal{L}^{\dagger}\varphi\bigr)\,d\tau
                  =\int_{G}\bigl(\psi_{\epsilon}*u\bigr)\,\mathcal{L}^{\dagger}\varphi\,d\tau
                  =\int_{G}\bigl(\mathcal{L}u_{\epsilon}\bigr)\,\varphi\,d\tau,
              \end{equation*}
              where the last equality is integration by parts against the smooth bounded function
              $u_{\epsilon}=\psi_{\epsilon}*u$. Since $\varphi\in \mathcal{D}$ was arbitrary and
              $\mathcal{L}u_{\epsilon}$ is continuous, $\mathcal{L}u_{\epsilon}=0$ pointwise.
              The other direction follows from the direct consequence of the fact that $u_{\epsilon}\to u$
              weak-$*$ in $L^\infty(G)$ as $\epsilon\to 0$, and $\mathcal{L}^{\dagger}\psi\in L^1(G)$ for every
              $\psi\in \mathcal{D}(G)$,
              \begin{equation*}
                  \int_{G} u\,\mathcal{L}^{\dagger}\psi\,d\tau=
                  \lim_{\epsilon\to 0}\int_{G}u_{\epsilon}\,\mathcal{L}^{\dagger}\psi \, d\tau
                  =\lim_{\epsilon\to 0}\langle \mathcal{L}u_{\epsilon},\psi\rangle=0.
              \end{equation*}

        \item[2.]
              Assume $R_{s}u=u$ a.e. By \eqref{eq:conv-adjoint} and
              property~(iii) above,
              \begin{equation*}
                  \int_{G}(R_{s}u_{\epsilon}-u_{\epsilon})\,\psi\,dh
                  =\int_{G}\bigl(\psi_{\epsilon}*(R_{s}u-u)\bigr)\,\psi\,dh
                  =\int_{G}(R_{s}u-u)\,(\tilde\psi_{\epsilon}*\psi)\,dh
                  =0
              \end{equation*}
              for all $\psi\in C_c^\infty(G)$, since $R_{s}u-u=0$ a.e., continuity of
              $R_{s}u_{\epsilon}-u_{\epsilon}$ gives $R_{s}u_{\epsilon}=u_{\epsilon}$ pointwise.

              The other direction follows similarly as in the proof of $1$.
    \end{itemize}
\end{proof}

\section{Proof via a Choquet--Deny formulation}\label{sec:choquet}

We explain our approach for proving Theorem~\ref{thm:main} via the Choquet--Deny formulation. Note that this proof allows us to cover the nilpotent case as well as the compact case. We state a result of B\'enard and Frisch~\cite{BenardFrisch} that is critical to our proof of the nilpotent case.

\begin{theorem}[B\'enard--Frisch~\cite{BenardFrisch}]\label{thm:benard-frisch}
    Let $G$ be a connected nilpotent Lie group, and suppose its iterated commutator subgroup
    $[G,[G,G]]$ is relatively compact. If $\nu$ is a probability measure on $G$ whose support
    generates a dense subgroup, then every bounded measurable $h\colon G\to\mathbb{R}$ satisfying
    \begin{equation*}
        h(x)=\int_G h(x\star y)\,\nu(dy)\qquad\text{for all }x\in G,
    \end{equation*}
    is necessarily constant.
\end{theorem}

\begin{remark}
    B\'enard and Frisch's convolution equation is the classical Choquet--Deny equation after
    reflecting the measure. Starting from
    \begin{equation*}
        h(x)=\int_G h(x\star y)\,\nu(dy)\qquad\text{for $x\in G$},
    \end{equation*}
    if we take $\check\nu$ to be the pushforward of $\nu$ under inversion $y\mapsto y^{-1}$
    \begin{equation*}
        \check\nu(A)=\nu(A^{-1}), \quad \text{for Borel $A$,}
    \end{equation*}
    so that
    \begin{equation*}
        \int_G\varphi(z)\,\check\nu(dz)=\int_G\varphi(y^{-1})\,\nu(dy),
    \end{equation*}
    for every test function $\varphi$. Taking $\varphi(z)=h(x\star z^{-1})$ gives
    \begin{equation*}
        h(x)=\int_G h(x\star y)\,\nu(dy)=\int_G h(x\star z^{-1})\,\check\nu(dz)=(h*\check\nu)(x),
    \end{equation*}
    with the convention $(h*\sigma)(x):=\int_G h(x\star z^{-1})\,\sigma(dz)$. Hence
    $h=h*\check\nu$, the original Choquet--Deny equation.

    Since
    \begin{equation*}
        \supp\check\nu=(\supp\nu)^{-1},
    \end{equation*}
    generates the same closed subgroup as $\supp\nu$, the measure
    $\check\nu$ is adapted, i.e.\ $\overline{\langle \supp \check\nu\rangle}=G$, if and only if $\nu$ is. Hence $h$ is constant when either variant of the Choquet--Deny equation is satisfied.
\end{remark}

We now outline our approach to proving Theorem~\ref{thm:main}. We assume that our solution $u$ is continuous thanks to Proposition~\ref{prop:smooth_soln}. Note that Corollary~\ref{cor:choquet-deny} enables us to consider the following problem
\begin{equation*}
    u=u*\rho_{t}, \quad \text{for $t\geq 0$,}
\end{equation*}

We now explain our approach to proving Theorem \ref{thm:main}.

\begin{itemize}
    \item[\emph{Step 1.}] We show that as a consequence of~\eqref{eqn:convolution_eqn1}, we have that $u=u*\rho$.
\end{itemize}
We aim to show that $u(g\star p)=u(g)$ for $g\in G$ and $p\in \overline{\langle \supp \rho \rangle}$. We treat the compact and nilpotent cases separately.
\begin{itemize}
    \item[\emph{Step 2.}] We show this under the assumption that $G$ is compact via the It\^o--Kawada theorem, see~\cite{MR3462} or~\cite{MR4330970}.
    \item[\emph{Step 3.}] We directly apply the main result of B\'enard and Frisch~\cite{BenardFrisch} to consider the case when $G$ is a nilpotent Lie group with the property that $[G,[G,G]]$ is relatively compact.
\end{itemize}

We now denote $P=\overline{\langle\supp\rho\rangle}$. The goal is now to show that  $P=\overline{\langle\exp\mathfrak h_0\cup\supp\mu\rangle}$,
$\mathfrak h_0=\Lie\langle b,\operatorname{ran}(a)\rangle$. By~\cite[Proposition 6.7]{MR2060091},
$\supp\rho=\overline{\bigcup_{t>0}\supp\rho_t}$, so $P$ is the smallest closed subgroup $H$
with $\rho_t(H)=1$ for all $t$. We therefore ask which closed subgroups carry the whole
semigroup. We fix a closed subgroup $H$, choose a basis $X_1,\dots,X_d$ of $\Lie(G)$ with
$X_1,\dots,X_m$ spanning $\Lie(H)$ and $X_{m+1},\dots,X_d$ transverse, and write
$x_i(g)$ for the $i$-th coordinate of $\log g$, so that near $e$ the subgroup $H$ is
$\{x_{m+1}=\dots=x_d=0\}$.

\begin{itemize}
    \item[\emph{Step 4.}] We prove that $\rho_t(H)=1$ for all $t$ forces
          \begin{equation}\label{eqn:triple_assumptions}
              \supp\mu\subseteq H, \quad \operatorname{ran}(a)\subseteq\Lie(H), \quad  b\in\Lie(H).
          \end{equation}

    \item[\emph{Step 5.}] Conversely, we show that if the above 3 conditions hold, the characteristics
          restrict to a L\'evy triplet on $H$ whose generator coincides with $\mathcal{L}$ on
          $C_c^\infty(G)$. Hunt's uniqueness~\cite{MR79232} then forces $\rho_t(H)=1$ for all $t$.

    \item[\emph{Step 6.}] We prove $P=\overline{\langle\exp\mathfrak h_0\cup\supp\mu\rangle}$,
          $\mathfrak h_0=\Lie\langle b,\operatorname{ran}(a)\rangle$ which completes our proof.
\end{itemize}

\begin{proof}[Proof of Theorem~\ref{thm:main} via a Choquet--Deny type approach]
We are now ready to complete the proof of Theorem~\ref{thm:main}.

\step{Step 1.}
Integrating~\eqref{eqn:convolution_eqn1} with respect to $e^{-t}$, we obtain
\begin{equation*}
    \begin{split}
        u(g) & =\int_0^\infty e^{-t}u(g)\,dt                                        \\
             & =\int_0^\infty e^{-t}\,(u*\rho_t)(g)\,dt                             \\
             & =\int_0^\infty e^{-t}\!\left(\int_G u(g\star h)\,\rho_t(dh)\right)dt \\
             & =\int_G u(g\star h)\,\rho(dh)=(u*\rho)(g),\qquad g\in G,
    \end{split}
\end{equation*}
which is our Choquet--Deny equation
\begin{equation}\label{eqn:choquet_deny}
    u=u*\rho.
\end{equation}
Note, for the above steps to be valid, we used the Fubini theorem and the fact that $u$ is bounded.

\step{Step 2.}
We begin with the compact setting which uses the It\^o--Kawada theorem, see~\cite{MR3462} or~\cite{MR4330970}. We present a proof using a formulation given in ~\cite{MR3462}. This is standard for readers familiar with Choquet--Deny theory, but we include it for clarity.

Using the fact that $u=u*\rho=\cdots=u*\rho^{*k}$ for $k\in\mathbb{N}$, we have that
\begin{equation*}
    u=u*A_{N}, \qquad A_{N}:=\frac{1}{N}\sum_{k=0}^{N-1}\rho^{*k}\quad(N\in\mathbb{N}),
\end{equation*}
where $\rho^{*k}$ is the $k$-th convolution power of $\rho$. By an application of the It\^o--Kawada theorem we have that
\begin{equation*}
    u=u*A_{N}=\lim_{N\rightarrow \infty} u*A_{N}=u*m_{P},
\end{equation*}
where $m_{P}$ is the Haar measure supported on $\overline{\langle \supp \rho \rangle}$.

For each $g$, expanding the square and using $\int_P u(g\star p)\,dm_P(p)=u(g)$ and $\int_P dm_P=1$,
\begin{equation*}
    \int_P |u(g\star p)-u(g)|^2\,dm_P(p)=\int_P |u(g\star p)|^2\,dm_P(p)-u(g)^2 .
\end{equation*}
Integrate over $g\in G$ against $m_G$, by Fubini and the right-invariance of $m_G$,
$\int_G |u(g\star p)|^2\,dm_G(g)=\|u\|_2^2$ for each $p$, so
\begin{equation*}
    \int_G\!\int_P |u(g\star p)-u(g)|^2\,dm_P(p)\,dm_G(g)=\|u\|_2^2-\|u\|_2^2=0 .
\end{equation*}
The integrand is nonnegative, hence vanishes a.e., and we necessarily have that
\begin{equation*}
    u(g\star p)=u(g) \quad \text{for $m_G$-a.e.\ $g$ and $m_P$-a.e.\ $p$}.
\end{equation*}
Since $u$ is continuous, we have completed the step for when $G$ is compact.

\step{Step 3.} We now consider the nilpotent case under the assumption that $[G,[G,G]]$ is relatively compact.

Fix $g\in G$ and set
$u_g(p)=u(g\star p)$ for $p\in P$. Since $\rho$ is carried by $P$,
\eqref{eqn:choquet_deny} restricts to $u_g=u_g*\rho$ on $P$, so $u_g$ is a bounded
$\rho$-harmonic function on $P$, and $\overline{\langle \supp \rho \rangle }=P$ by construction. As a closed subgroup of the nilpotent Lie group $G$, the group $P$ is itself a nilpotent Lie group, and
its iterated commutator satisfies $[P,[P,P]]\subseteq[G,[G,G]]$ and is therefore relatively compact. Hence
$P$ is a Choquet--Deny group and $u_g$ is constant on $P$. When $P$ is connected this is exactly
Theorem~\ref{thm:benard-frisch}. In general the identity component $P^{\circ}$ is a connected nilpotent Lie
group with relatively compact iterated commutator, hence Choquet--Deny by Theorem~\ref{thm:benard-frisch},
while the component group $P/P^{\circ}$ is finitely generated nilpotent and hence Choquet--Deny by the
theorem of Dynkin and Malyutov~\cite{DynkinMalyutov} (see also~\cite{FHTV}). Since $P$ is nilpotent, these
two facts combine to give the Choquet--Deny property for $P$. Thus $u(g\star p)=u_g(e)=u(g)$ for all $p\in P$, and as
$g$ was arbitrary, $u$ is right-invariant under $P$.

\step{Step 4.}
Suppose $\rho_t(H)=1$ for all $t$. Any $f\in C_c^\infty(G)$ vanishing on $H$ then has
$\int_G f\,d\rho_t=0$ for all $t$ and $f(e)=0$, so
\begin{equation*}
    \mathcal{L}f(e)=\lim_{t\to0}\frac1t\Big(\int_G f\,d\rho_t-f(e)\Big)=0 .
\end{equation*}
We verify that the sets are contained in the support by choosing $f$ carefully.

\medskip\noindent\textit{(a) $\supp\mu\subseteq H$}\quad\ \\
Take any $f\ge0$ with $\supp f\subseteq G\setminus H$. As $e\in H$, the point $e$ is outside
the closed set $\supp f$, so $f\equiv0$ on a neighborhood of $e$. In turn, we have that
\begin{equation*}
    f(e)=0,\quad (X_kf)(e)=0,\quad  \text{and}\quad (X_pX_qf)(e)=0.
\end{equation*}
It follows that
\begin{equation*}
    \mathcal{L}f(e)=\int_G f(\tau)\,\mu(d\tau)=0.
\end{equation*}
This holds for every nonnegative $f$ supported in $G\setminus H$, which forces $\mu(G\setminus H)=0$, i.e.
$\supp\mu\subseteq H$.

\medskip\noindent\textit{(b) $\operatorname{ran}(a)\subseteq\Lie(H)$}\quad\ \\
Fix transverse $i,j$ and take $f=x_ix_j\chi$. Using the Leibniz
rule and $x_i(e)=x_j(e)=0$, we have that
\begin{equation*}
    f(e)=0,\qquad (X_kf)(e)=(X_kx_i)(e)\,x_j(e)+x_i(e)\,(X_kx_j)(e)=0,
\end{equation*}
\begin{equation*}
    (X_pX_qf)(e)=(X_px_i)(e)(X_qx_j)(e)+(X_qx_i)(e)(X_px_j)(e)
    =\delta_{pi}\delta_{qj}+\delta_{pj}\delta_{qi}.
\end{equation*}
Hence, by direct computation, we have that
\begin{align*}
    \mathcal{L}f(e) & =\sum_{p,q}a_{pq}\big(\delta_{pi}\delta_{qj}+\delta_{pj}\delta_{qi}\big)
    +\int_G f\,d\mu,                                                                          \\
                   & =(a_{ij}+a_{ji})+\int_G f\,d\mu,                                         \\
                   & =2a_{ij}+\int_G f\,d\mu,                                                 \\
                   & =2a_{ij},
\end{align*}
since, from (a), we have $\supp\mu\subseteq H$ and $f$ vanishes on $H$.

Since $\mathcal{L}f(e)=0$ we necessarily have that $a_{ij}=0$ for every pair of transverse indices. Taking $j=i$ shows
$a_{ii}=0$ for every transverse $i$. As $a$ is symmetric and positive semidefinite, each
$2\times2$ principal minor is nonnegative, so $a_{ii}a_{kk}-a_{ik}^2\ge0$. Hence $a_{ik}=0$ for all $k$. Thus every transverse row of $a$ vanishes, which is
precisely $\operatorname{ran}(a)\subseteq\Lie(H)$.

\medskip\noindent\textit{(c) $b\in\Lie(H)$.}\quad\ \\
Fix a transverse index $i$ and set $f=x_i\chi$. Since $\chi\equiv1$ near $e$, we have $f=x_i$
on a neighborhood of $e$, so $f(e)=0$ and the jet of $f$ at $e$ coincides with that of $x_i$
\begin{equation*}
    (X_kf)(e)=(X_kx_i)(e)=\delta_{ki},\qquad
    (X_pX_qf)(e)=(X_pX_qx_i)(e)=\tfrac12\gamma^i_{pq},
\end{equation*}
where $\gamma^i_{pq}$ is determined by the relation $[X_p,X_q]=\sum_r\gamma^r_{pq}X_r$, see also \cite[Proposition 3.7]{MR722297}.
By direct calculation,
\begin{equation*}
    \mathcal{L}f(e)=\sum_{p,q}a_{pq}(X_pX_qf)(e)+\sum_p b_p(X_pf)(e)
    +\int_G\Big(f(\tau)-f(e)-\sum_k x_k(\tau)(X_kf)(e)\Big)\mu(d\tau),
\end{equation*}
and evaluate its three parts.
The first part is given by
\begin{equation*}
    \sum_{p,q}a_{pq}(X_pX_qf)(e)=\sum_{p,q}a_{pq}\,\tfrac12\gamma^i_{pq}=0,
\end{equation*}
since $a_{pq}=a_{qp}$ while $\gamma^i_{pq}=-\gamma^i_{qp}$, so relabelling $p\leftrightarrow q$, the sum is equal to its own negative.

The drift term is given by
\begin{equation*}
    \sum_p b_p(X_pf)(e)=\sum_{p}b_p \delta_{pi}=b_i.
\end{equation*}
The jump term is also given by
\begin{align*}
    \int_G\Big(f(\tau)-f(e)-\sum_k x_k(\tau)(X_kf)(e)\Big)\mu(d\tau)
     & =\int_G\Big(x_i(\tau)\chi(\tau)-0-\sum_k x_k(\tau)\delta_{ki}\Big)\mu(d\tau) \\
     & =\int_G\big(x_i(\tau)\chi(\tau)-x_i(\tau)\big)\mu(d\tau)                     \\
     & =\int_G x_i(\tau)\big(\chi(\tau)-1\big)\mu(d\tau)                            \\
     & =0 .
\end{align*}
Here we used $f(\tau)=x_i(\tau)\chi(\tau)$ and $f(e)=x_i(e)\chi(e)=0$ (as $x_i(e)=0$), and
the compensator $\sum_k x_k(\tau)(X_kf)(e)=\sum_k x_k(\tau)\delta_{ki}=x_i(\tau)$. The final
integral vanishes because, by part (a), $\supp\mu\subseteq H$, and the transverse coordinate
$x_i$ vanishes on $H$, hence $x_i=0$ $\mu$-almost everywhere, so the integrand is $0$
$\mu$-a.e.

Since $\mathcal{L}f(e)=0$ we get
$b_i=0$. As $i$ ranges over the transverse indices, $b\in\Lie(H)$.

\step{Step 5.}
Conversely, assume $\supp\mu\subseteq H$, $\operatorname{ran}(a)\subseteq\Lie(H)$ and
$b\in\Lie(H)$, and keep the adapted basis with $X_1,\dots,X_m$ spanning $\Lie(H)$.

The three assumptions, given in~\eqref{eqn:triple_assumptions}, tell us when the characteristics have no
transverse part.

Note that $a_{jk}=0$ whenever $j>m$ or $k>m$, $b_i=0$ for $i>m$, and $\mu$ is carried
by $H$. Hence the restricted data
\begin{equation*}
    \big(a|_{\Lie(H)},\,b|_{\Lie(H)},\,\mu|_H\big),\qquad
    a|_{\Lie(H)}=(a_{jk})_{j,k\le m},\quad b|_{\Lie(H)}=(b_i)_{i\le m},
\end{equation*}
is a L\'evy triplet on $H$ since
\begin{itemize}
    \item $a|_{\Lie(H)}$ is symmetric positive semidefinite,
    \item $\mu|_H$
          is a L\'evy measure on $H$, since for any neighborhood $V$ of $e$ one has
          $\int_V\sum_{i\le m}x_i^2\,d\mu\le\int_V\sum_i x_i^2\,d\mu<\infty$ and $\mu(V^c)<\infty$.
\end{itemize}
In turn, by Hunt's theorem~\cite{MR79232} this triplet generates a unique convolution semigroup
$(\rho^H_t)_{t\ge0}$ on $H$, with generator $\widetilde{\mathcal{L}}$.

For $f\in C_c^\infty(G)$, subtracting the two generators leaves only the transverse part of
each characteristic
\begin{equation*}
    \mathcal{L}f-\widetilde{\mathcal{L}}f
    =\sum_{j>m\ \text{or}\ k>m}a_{jk}X_jX_kf
    +\sum_{i>m}b_i\,X_if
    -\int_G\sum_{i>m}x_i(\tau)\,(X_if)\,\mu(d\tau).
\end{equation*}
Observe that
\begin{equation*}
    \sum_{j>m\ \text{or}\ k>m}a_{jk}X_jX_kf=0,
\end{equation*}
since $a_{jk}=0$ once $j$ or $k$ is transverse
($\operatorname{ran}(a)\subseteq\Lie(H)$).

Moreover, we have that
\begin{equation*}
    \sum_{i>m}b_i\,X_if=0,
\end{equation*}
since
\begin{equation*}
    b_i=0,\quad \text{for $i>m$, recall that $b\in\Lie(H)$},
\end{equation*}
and
\begin{equation*}
    \int_G\sum_{i>m}x_i(\tau)\,(X_if)\,\mu(d\tau)=0,
\end{equation*}
since $x_i=0$ on $\supp\mu\subseteq H$ for transverse $i$,
so the integrand is $0$ $\mu$-a.e.

Thus $\mathcal{L}=\widetilde{\mathcal{L}}$ on $C_c^\infty(G)$.

As a convolution semigroup is determined by the restriction of its generator to $C_c^\infty$
(uniqueness in Hunt's theorem~\cite{MR79232}), $\rho_t=\rho^H_t$ for every $t$. Moreover, since each
$\rho^H_t$ is carried by $H$, we conclude $\rho_t(H)=1$ for all $t$.

\step{Step 6.}
By Steps 4 and 5, for a closed subgroup $H$ one has $\rho_t(H)=1$
for all $t$ if and only if
\begin{equation*}
    \supp\mu\subseteq H,\qquad \operatorname{ran}(a)\subseteq\Lie(H),\qquad b\in\Lie(H) .
\end{equation*}
Recalling from Step~3 that $P=\overline{\langle\supp\rho\rangle}$ is the smallest closed
subgroup with $\rho_t(P)=1$ for all $t$, it follows that $P$ is the smallest closed subgroup
satisfying these three conditions. We identify it explicitly. Set
$\mathfrak h_0:=\Lie\langle b,\operatorname{ran}(a)\rangle$, the Lie subalgebra generated by the
drift and the range of the diffusion.

If a closed subgroup $H$ satisfies the three conditions, then $\operatorname{ran}(a)\subseteq\Lie(H)$
and $b\in\Lie(H)$, since $\Lie(H)$ is a subalgebra, $\mathfrak h_0\subseteq\Lie(H)$,
whence $\overline{\langle\exp\mathfrak h_0\rangle}\subseteq H$, and also $\supp\mu\subseteq H$.
Therefore
\begin{equation*}
    \overline{\langle\exp\mathfrak h_0\cup\supp\mu\rangle}\subseteq H .
\end{equation*}
Conversely, $\overline{\langle\exp\mathfrak h_0\cup\supp\mu\rangle}$ itself satisfies the three
conditions, in that it contains $\supp\mu$, its Lie algebra contains $\mathfrak h_0\supseteq\operatorname{ran}(a)$,
and it contains $b$. Being the smallest such subgroup, it equals $P$
\begin{equation*}
    P=\overline{\langle\supp\rho\rangle}
    =\overline{\big\langle\,\exp\mathfrak h_0\cup\supp\mu\,\big\rangle},\qquad
    \mathfrak h_0=\Lie\langle b,\operatorname{ran}(a)\rangle .
\end{equation*}

Finally, suppose $P=G$, that is, $\exp\mathfrak h_0$ together with $\supp\mu$
generates a dense subgroup. Combined with the right-invariance of Step~2 (compact case) or
Step~3 (nilpotent case), $u(g\star p)=u(g)$ for all $g\in G$ and $p\in P=G$.
Taking $g=e$ gives $u(p)=u(e)$ for every $p\in G$, i.e. $u$ is constant.

Conversely, suppose $P\neq G$. Since $G$ is connected, a proper closed subgroup has positive
codimension, so the left coset space $G/P$ is a smooth manifold of dimension at least one. Choose any
nonconstant bounded smooth function $\varphi$ on $G/P$ and set $u(g)=\varphi(gP)$. Since $\rho_t$ is carried by
$P$, every $h\in\supp\rho_t$ lies in $P$, so $g\star h$ and $g$ determine the same left coset and
$u(g\star h)=u(g)$. Hence $u=u*\rho_t$ for all $t\ge0$, and differentiating at $t=0$ gives $\mathcal{L}u=0$.
Thus $u$ is a nonconstant bounded solution, so not every bounded solution is constant when $P\neq G$.
\end{proof}

\section{A special case of Theorem~\ref{thm:main} via propagation of maxima}
\label{sec:propagation}

In this section, we aim to prove a special case of Theorem~\ref{thm:main} via a propagation of maxima approach due to Alibaud, del~Teso, Endal and Jakobsen~\cite{MR4149690}. We aim to show why we cannot get the full generality of Theorem~\ref{thm:main} via this approach. Namely, we can only prove the following weaker version of Theorem~\ref{thm:main}.

\begin{theorem}\label{thm:compact-propagation}
    Let $G$ be a connected compact Lie group and $\mathcal{L}$ a left-invariant L\'evy operator
    with characteristics $(a,b,\mu)$. Set
    $\mathfrak h_0=\Lie\langle b,\operatorname{ran}(a)\rangle$ and
    $P=\overline{\langle\exp\mathfrak h_0\cup\supp\mu\rangle}$. Then every bounded distributional
    solution $u$ of $\mathcal{L}u=0$ is constant if and only if $P=G$.
\end{theorem}

Note that we do not claim that $u$ is right-invariant when $P \ne G$ unlike in Theorem~\ref{thm:main}. We present our proof and highlight the limitations that we face in obtaining Theorem~\ref{thm:main}.

By Proposition~\ref{prop:smooth_soln} we may assume $u\in C^\infty(G)$ solves $\mathcal{L}u=0$
pointwise. Since $G$ is compact and $u$ is continuous, $u$ attains its maximum $M=\max_G u$. Set
$F=\{g\in G:u(g)=M\}$, a nonempty closed set.

The argument runs in four steps.
\begin{itemize}
    \item[\emph{Step 1.}] We aim to show that, for a given $\overline{g}\in F$, we have that
          \begin{equation} \label{eqn:jump_part}
              \int_G\big(u(\bar g\star\tau)-M\big)\mu(d\tau)=0,
          \end{equation}
          and
          \begin{equation}\label{eqn:diffusion_part}
              \sum_{i,j}a_{ij}X_iX_ju(\bar g)+\sum_i b_iX_iu(\bar g)=0.
          \end{equation}
    \item[\emph{Step 2.}] As a consequence of~\eqref{eqn:jump_part} we show that $\overline{g}\star \supp \mu \in F$.
    \item[\emph{Step 3.}] On its maximum set $u$ is
          a subsolution of the local part of $\mathcal{L}$. Hence, by Bony's strong maximum principle~\cite{Bony},
          together with the compactness of $G$, $F$ is invariant under $\exp\mathfrak h_0$.
    \item[\emph{Step 4.}]  Combining Steps 2 and 3, $F$ contains the right coset $\bar g\star P$. Hence,
          $u$ is constant if and only if $P=G$.
\end{itemize}

\begin{proof}[Proof of Theorem~\ref{thm:compact-propagation}]
    \step{Step 1.}
    Fix $\bar g\in F$. At this maximum $X_iu(\bar g)=0$ and the symmetric Hessian is negative
    semidefinite, so in $\mathcal{L}u(\bar g)=0$ the drift and compensator terms drop out, leaving
    \begin{equation*}
        0=\sum_{i,j}a_{ij}X_iX_ju(\bar g)+\int_G\big(u(\bar g\star\tau)-M\big)\mu(d\tau).
    \end{equation*}
    The first term is $\le0$ ($a$ is positive semidefinite and the Hessian negative semidefinite) and
    the second is $\le0$ (as $u\le M$), so both vanish. This is~\eqref{eqn:jump_part}, and, since
    $\sum_i b_iX_iu(\bar g)=0$, we also obtain ~\eqref{eqn:diffusion_part}.

    \step{Step 2.}
    In~\eqref{eqn:jump_part} the integrand $u(\bar g\star\tau)-M$ is $\le0$, so it vanishes for
    $\mu$-a.e.\ $\tau$, hence for all $\tau\in\supp\mu$ by continuity. Thus $\bar g\star\supp\mu\subseteq F$.

    \step{Step 3.}
    Every point of $F$ is a maximum, so~\eqref{eqn:diffusion_part} gives $\mathcal{L}^{\mathrm{loc}}u=0$ on
    $F$, where $\mathcal{L}^{\mathrm{loc}}=\sum_{i,j}a_{ij}X_iX_j+\sum_i b_iX_i$. In particular, $u$ is a
    subsolution of $\mathcal{L}^{\mathrm{loc}}$ attaining its maximum on $F$. By Bony's strong maximum principle~\cite{Bony},
    $F$ is invariant under the diffusion flows $\exp(tY)$, for $Y\in\operatorname{ran}(a)$, $t\in\R$, and the forward drift flow $\exp(tb)$ ($t\ge0$). Since $G$ is
    compact, $\{\exp(tb):t\ge0\}$ is dense in the torus $\overline{\{\exp(tb):t\in\R\}}$, so the drift
    flow is two-sided. Closing under brackets, $\bar g\star\exp\mathfrak h_0\subseteq F$ with
    $\mathfrak h_0=\Lie\langle b,\operatorname{ran}(a)\rangle$.
    \step{Step 4.}
    Being closed and right-invariant under $\supp\mu$ and $\exp\mathfrak h_0$, $F$ is invariant under the
    subgroup they generate, so $\bar g\star P\subseteq F$ with
    $P=\overline{\langle\exp\mathfrak h_0\cup\supp\mu\rangle}$. If $P=G$ then $F=G$ and $u$ is constant.
    If $P\ne G$ there are nonconstant bounded solutions, constructed as in the proof of Theorem~\ref{thm:main} (for instance Example~\ref{ex:su2}). Hence $u$ is
    constant if and only if $P=G$.
\end{proof}

\subsection*{A comparison to the original approach of Alibaud, del~Teso, Endal and Jakobsen}
The approach of Alibaud, del~Teso, Endal and Jakobsen~\cite{MR4149690} establishes more than the
dichotomy of Theorem~\ref{thm:compact-propagation}. It shows that every bounded solution is
right-invariant under $P$, not merely constant when $P=G$. We are not able to carry this stronger
conclusion over to a general compact or nilpotent Lie group, and we explain why.

On $\R^d$ their argument runs as follows. By Proposition~\ref{prop:smooth_soln} take
$u\in C^\infty_b(\R^d)$ solving $\mathcal{L}u=0$ pointwise, and for a jump direction $\bar z\in\supp\mu$
form the increment $w(x)=u(x+\bar z)-u(x)$. Since translation commutes with $\mathcal{L}$, $w$ is again a
bounded solution. Choosing $x_n$ with $w(x_n)\to\sup w$, each $u_n:=u(\cdot+x_n)$ is again a
bounded solution, and the mollified bounds $\|D^k u_n\|_\infty\le\|u\|_\infty\|D^k\rho_\varepsilon\|_{L^1}$
are uniform in $n$ because, on $\R^d$, translation commutes with differentiation. By Arzel\`a--Ascoli a
subsequence converges in $C^2_{\mathrm{loc}}$ to a solution $u_\infty$, and
$w_\infty=u_\infty(\cdot+\bar z)-u_\infty$ attains $\sup w$ at the origin, so a value approached but never attained by the original increment becomes a genuine interior maximum of $w_\infty$. Propagation of the maximum there, with a
telescoping argument, forces $\sup w=0$, that is $u(\cdot+\bar z)=u$. The same applies in the diffusion and
drift directions, giving the full right-invariance.

Two steps of this argument fail on a general group. First, it is essential that the increment
$w=u(\cdot+\bar z)-u$ again solves the equation. On $\R^d$ this is translation invariance, but on a
non-commutative group the right translate $R_{\bar z}u$ solves the operator conjugated by
$\operatorname{Ad}(\bar z^{-1})$ rather than $\mathcal{L}$ itself, so $w$ need not be a solution and the
increment method breaks down. This is why, even on a compact group, propagation reaches only the maximum
set and the dichotomy of Theorem~\ref{thm:compact-propagation}, not the right-invariance. Second,
manufacturing a maximizer rests on the basepoint-uniform derivative bounds above, which can fail on a
nilpotent group, as the following example shows.

\begin{example}\label{ex:heisenberg}
    Take the Heisenberg group $\mathbb{H}^1=\mathbb{R}^3$ with the group operation given by
    \begin{equation*}
        (x,y,t)\star(x',y',t')=(x+x',y+y',t+\tfrac12(xy'-yx')),
    \end{equation*}
    with the left-invariant fields
    \begin{equation*}
        X=\partial_x-\tfrac y2\partial_t, \qquad Y=\partial_y+\tfrac x2\partial_t.
    \end{equation*}

    We take $\mathcal{L}=X^2$ and the function $u$
    \begin{equation*}
        u(x,y,t)=\sin(y^2), \quad \text{for $(x,y,t)\in \mathbb{H}^1$.}
    \end{equation*}
    Since $u$ depends only on $y$, $Xu=0$, so $u$ is a bounded solution of
    $\mathcal{L}u=0$. It is nonconstant, and right-invariant under the subgroup
    $P=\{(s,0,0)\}$ since
    \begin{equation*}
        u\left( (x,y,t)\star(s,0,0)\right )=\sin(y^2)=u(x,y,t).
    \end{equation*}
    Run the manufacturing step above, sliding along $g_n=(0,y_n,0)$ with
    $y_n\to\infty$ and consider the translates $u\circ L_{g_n}$. The transverse
    left-invariant derivative is
    \begin{equation*}
        Yu=2y\cos(y^2),\qquad \|Yu\|_{L^\infty(\mathbb{H}^1)}=\infty.
    \end{equation*}
    In turn, we cannot conclude that $\{u\circ L_{g_n}\}$ is equicontinuous, nor extract a $C^2_{\mathrm{loc}}$-convergent subsequence, as Arzel\`a--Ascoli fails, and no maximizer can be manufactured.
\end{example}

\section{A harmonic analysis approach of Theorem~\ref{thm:main} in the case of compact Lie groups}\label{sec:harmonic}
\label{sec:fourier}

In this section, we provide a harmonic analysis approach to Theorem~\ref{thm:main}, in the case where $G$ is compact.

In this setting we follow the Fourier-analytic approach to L\'evy processes on compact Lie groups developed by Liao~\cite{MR2060309} and Applebaum~\cite[Chapter 5]{MR3243650} (see also~\cite{MR3335018}), in which the generator is realized as a Fourier multiplier acting on the matrix coefficients of the irreducible representations. We use the Peter--Weyl theory of harmonic analysis on compact Lie groups~\cite{MR1512386}, cf. the books of Hall~\cite[Chapter 12]{MR3331229} and Sepanski~\cite[\S 3.3]{Sepanski}.

The Peter--Weyl theory decomposes the $L^2(G)$ space into isotypic components
\begin{equation*}
    L^2(G) = \bigoplus_{\pi \in \hat{G}}V_\pi.
\end{equation*}
The Fourier transform, for $f\in \mathcal{D}(G)$, is understood as
\begin{equation*}
    \mathcal{F}(f)(\pi) = \int_G f(g)\pi(g^{-1}) dg, \quad \text{for $\pi \in \hat{G}$,}
\end{equation*}
where $dg$ is the Haar measure on $G$. We write $\mathcal{S}(\hat{G})$ for the \emph{Schwartz space} on the unitary dual $\hat{G}$, consisting of the matrix fields $h=(h(\pi))_{\pi\in\hat{G}}$, with $h(\pi)\in\mathbb{C}^{\dim(\pi)\times\dim(\pi)}$, that decay rapidly in the sense that $\sup_{\pi\in\hat{G}}(1+\lambda_\pi)^{N}\lVert h(\pi)\rVert<\infty$ for every $N\in\mathbb{N}$, where $\lambda_\pi\geq0$ is the eigenvalue by which the Laplace--Beltrami operator acts on the matrix coefficients of $\pi$; the Fourier transform is a linear bijection from $\mathcal{D}(G)$ onto $\mathcal{S}(\hat{G})$. The Fourier inverse is given as
\begin{equation*}
    \mathcal{F}^{-1}(h)(g) = \sum_{\pi \in \hat{G}} \dim(\pi) \operatorname{Tr}(\pi(g)h(\pi)),
    \quad \text{for $g \in G$, $h \in \mathcal{S}(\hat{G}).$}
\end{equation*}

In the compact case the components $V_\pi$ are all finite dimensional vector spaces, thus the representations are finite dimensional homomorphisms from $G \to \operatorname{GL}(V_\pi),$ ($\operatorname{GL}(V_\pi)$ is the space of invertible matrices acting on $V_\pi$).

The corresponding distributional definition of the Fourier transform is given as
\begin{equation*}
    \langle \mathcal{F}u, \varphi \rangle_{\hat{G}} = \langle u, (\iota \circ \mathcal{F}^{-1})(\varphi) \rangle_{G},
\end{equation*}
where $\varphi \in \mathcal{S}(\hat{G})$ and $(\iota\circ \varphi)(x) = \varphi(x^{-1})$ for $\varphi \in \mathcal{D}(G)$.

The pairing in $\hat{G}$ is given as
\begin{align*}
    \langle u, \varphi \rangle_{\hat{G}} = \sum_{\pi \in \hat{G}} \dim(\pi) \operatorname{Tr}(u(\pi)\varphi(\pi)),
\end{align*}
where $\operatorname{Tr}$ is the usual trace operator for matrices.

From \cite[Theorem 5.5.1]{MR3243650} (cf.~\cite{MR2060309}), we have that
\begin{equation*}
    \widehat{\mathcal{L}f}(\pi)=\hat{\mathcal{L}}(\pi)\hat{f}(\pi), \quad \text{for $\pi \in \hat{G},$}
\end{equation*}
where
\begin{equation*}
    \begin{split}
        \hat{\mathcal{L}}(\pi) = & \sum_{i,j=1}^d a_{ij}\, d\pi(X_i)\, d\pi(X_j)+ \sum_{i=1}^d b_i\, d\pi(X_i) + \int_{G \setminus \{e\}} \left( \pi(\tau) - I - \sum_{i=1}^d x_i(\tau)\, d\pi(X_i) \right) \mu(d\tau), \\
    \end{split}
\end{equation*}
and $d\pi(X_i)$ is derived representation of $\pi(X_i)$.

We now outline our approach proving Theorem~\ref{thm:main}, in the compact setting.

\begin{itemize}
    \item[\emph{Step 1.}] We aim to show, using the notion of solution given in Definition~\ref{def:notion_of_solution}, that we have
          \begin{equation*}
              \hat{\mathcal{L}}(\pi)\hat{u}(\pi)=0, \quad \text{for all $\pi \in \hat{G}$,}
          \end{equation*}
          where $\hat{G}$ is the unitary dual of $G.$

    \item[\emph{Step 2.}] In this step, our goal is to characterize the kernel of $\hat{\mathcal{L}}(\pi)$. Firstly, for a given $v\in \ker \hat{\mathcal{L}}(\pi)$ (the kernel of $\hat{\mathcal{L}} (\pi)$) we necessarily have that \begin{equation*}\operatorname{Re}\langle \hat{\mathcal{L}}(\pi)v,v \rangle =0.\end{equation*} As a consequence, we determine that $v \in \ker \hat{\mathcal{L}}(\pi)$ must satisfy
          \begin{itemize}
              \item $d\pi(Y_l)v = 0$, for $Y_l = \sum_j\sigma_{jl}X_j\in\Lie(G)$,
              \item $(\pi(\tau) - I)v = 0$ for $\tau \in \supp(\mu)$,
          \end{itemize}
          More specifically, we have that $v$ must be invariant under the action of the subgroup $H_0$ generated by $\{\exp(tY_l) |\forall l,  t \in \mathbb{R} \}$ and $\supp(\mu)$.
    \item[\emph{Step 3.}] In this step we show that $v \in \ker \hat{\mathcal{L}}(\pi)$ must satisfy $d\pi(X_{\tilde{b}}) v = 0$ where
          \begin{equation*}
              X_{\tilde{b}} := \sum_{i=1}^d b_i X_i - \int_{G \setminus U} \left( \sum_{i=1}^d x_i(\tau) X_i \right) \mu(d\tau).
          \end{equation*}
          and $U$ is a neighborhood of the identity. Thus $v \in \ker(\hat{\mathcal{L}}(\pi))$ must be invariant under the subgroup $H$ which is generated by $H_0$ and $\{\exp(tX_{\tilde{b}}) | t \in \mathbb{R} \}$.

    \item[\emph{Step 4.}] In this step, we demonstrate that $X_{\tilde{b}}$ is in fact well defined  and show that the subgroup $H$ is not dependent on the choice of neighborhood $U$.
    \item[\emph{Step 5.}] We show $H$ periodicity of solutions $u$ and prove the converse, showing that any $H$ periodic distribution $u$ is a solution.
\end{itemize}

\begin{proof}[Proof of Theorem \ref{thm:main} in the compact case via Peter--Weyl]\
    \step{Step 1.}
    From the notion of solution, given in Definition~\ref{def:notion_of_solution}, we have that
    $\langle \mathcal{L}u, \varphi \rangle_G = \langle u, \mathcal{L}^\dagger\varphi \rangle_G$ for $\varphi \in \mathcal{D}(G)$. By then taking the distributional Fourier transform
    \begin{equation*}
        \langle u, \mathcal{L}^\dagger\varphi \rangle_G = 0 \iff \langle \hat{u}, \widehat{\iota \circ (\mathcal{L}^\dagger\varphi)}\rangle_{\hat{G}} = 0 .
    \end{equation*}

    The term
    \begin{align*}
        \widehat{\iota \circ (\mathcal{L}^\dagger\varphi)}(\pi) & = \int_G \big(\iota \circ (\mathcal{L}^\dagger\varphi)\big)(g)\, \pi(g^{-1})\, dg, \\
                                                                & = \int_G (\mathcal{L}^\dagger\varphi)(g^{-1})\, \pi(g^{-1})\, dg,                  \\
                                                                & = \int_G (\mathcal{L}^\dagger\varphi)(g)\, \pi(g)\, dg,                            \\
                                                                & = \int_G \varphi(g)\, \mathcal{L}(\pi)(g)\, dg,
    \end{align*}
    where we have used the inversion invariance of the Haar measure and the definition of the transpose.\\
    Computing the action of the operator on $\pi$ we see
    \begin{equation*}
        \mathcal{L}(\pi)(g) = \pi(g)\hat{\mathcal{L}}(\pi),\quad \text{for $(g,\pi)\in G\times \hat{G},$}
    \end{equation*}
    thus
    \begin{equation*}
        \int_G \varphi(g)\, \mathcal{L}\pi(g)\, dg = \left(\int_G \varphi(g)\, \pi(g)\, dg\right) \hat{\mathcal{L}}(\pi).
    \end{equation*}
    Moreover, we have that
    \begin{equation*}
        \langle \hat{u}, \widehat{\iota \circ (\mathcal{L}^\dagger\varphi)}\rangle_{\hat{G}} = \Big\langle \hat{u}, \Big(\int_G \varphi(g)\,\pi(g)\, dg\Big) \hat{\mathcal{L}}(\pi)\Big\rangle_{\hat{G}} = \langle \hat{u}, \widehat{(\iota \circ \varphi)}\, \hat{\mathcal{L}}\rangle_{\hat{G}}.
    \end{equation*}

    Evaluating the pairing we see
    \begin{align*}
        \langle \hat{u}, \widehat{(\iota \circ \varphi)}\, \hat{\mathcal{L}}\rangle_{\hat{G}} & = \sum_{\pi \in \hat{G}} \dim(\pi) \operatorname{Tr}\big( \hat{u}(\pi)\, \widehat{(\iota \circ \varphi)}(\pi)\, \hat{\mathcal{L}}(\pi)\big),      \\
                                                                                              & = \sum_{\pi \in \hat{G}} \dim(\pi) \operatorname{Tr}\big( \hat{\mathcal{L}}(\pi)\, \hat{u}(\pi)\, \widehat{(\iota \circ \varphi)}(\pi)\big) = 0 .
    \end{align*}
    As a consequence, we have that
    \begin{equation*}
        \hat{\mathcal{L}}(\pi)\hat{u}(\pi) = 0, \quad \text{for all $\pi \in \hat{G}.$}
    \end{equation*}
    Both these terms are matrices acting on the finite-dimensional representation space $V_\pi$, so this expression is $0$ when the image of $\hat{u}(\pi)$ is contained in the kernel of $\hat{\mathcal{L}}(\pi)$.
    \step{Step 2.} Our goal is to determine the kernel of $\hat{\mathcal{L}}(\pi).$\\
    Suppose $v \in V_\pi$ is in $\ker(\hat{\mathcal{L}}(\pi))$. Then necessarily, we have that
    \begin{equation*}
        \langle \hat{\mathcal{L}}(\pi)v,v\rangle_{V_\pi} = 0,
    \end{equation*}
    where $V_\pi$ is the vector space associated with the representation $\pi$.
    From the representation theory of compact Lie groups we know that the irreducible representations are finite dimensional. The operators $\hat{\mathcal{L}}(\pi)$ are thus matrices of dimension $\dim(\pi)$.

    First, we examine the real part. The term $\langle d\pi(X_i) v, v \rangle$ is purely imaginary as $d\pi(X_i)$ is skew-adjoint.
    As a consequence, we have that
    \begin{equation*}
        \operatorname{Re} \langle \hat{\mathcal{L}}(\pi) v, v \rangle =\sum_{j,k=1}^d a_{jk} \langle d\pi(X_j)\,d\pi(X_k) v, v \rangle + \int_{G \setminus \{e\}} \left( \operatorname{Re}\langle \pi(\tau) v, v \rangle - \langle v, v \rangle \right) \mu(d\tau).
    \end{equation*}

    The $d \times d$ matrix $(a_{jk})$ is nonnegative definite and symmetric, thus, along with the skew-adjointness of $d\pi(X_i)$, we have $\sum_{j,k=1}^d a_{jk} \langle d\pi(X_j)\,d\pi(X_k) v, v \rangle \leq 0$. Further, $(a_{jk})$ admits a decomposition $a = \sigma\sigma^{\mathsf T}$, i.e.\ writing out the matrix multiplication element-wise: $a_{jk} = \sum_l \sigma_{jl}\sigma_{kl}$. The second-order term becomes
    \begin{equation*}
        \sum_{j,k=1}^d a_{jk}\,d\pi(X_j)\,d\pi(X_k)=\sum_l\Big(\sum_j\sigma_{jl}\,d\pi(X_j)\Big)\Big(\sum_k\sigma_{kl}\,d\pi(X_k)\Big)=\sum_l d\pi(Y_l)^2,
    \end{equation*}
    where $Y_l:=\sum_j\sigma_{jl}X_j\in\Lie(G)$.
    The term
    \begin{equation*}
        \Big\langle \sum_l d\pi(Y_l)^2 v,v \Big\rangle = 0
    \end{equation*}
    when $d\pi(Y_l)v = 0$. By applying the exponential, we see that $\pi(\exp(tY_l))v = v$.
    The integrand in the second term is always negative as
    \begin{equation*}
        \operatorname{Re}\langle \pi(\tau)v, v \rangle - \lVert v\rVert^2 \le \lvert\langle \pi(\tau)v, v \rangle\rvert - \lVert v\rVert^2 \le \lVert\pi(\tau)v\rVert\,\lVert v\rVert - \lVert v\rVert^2\leq 0,
    \end{equation*}
    as a consequence of the Cauchy--Schwarz inequality and using the fact that $\pi$ is a unitary representation. Since $\pi(\tau)$ is unitary, $\lVert\pi(\tau)v\rVert = \lVert v\rVert$, thus the integrand is always $\leq 0$. Therefore as the integrand is negative the integral must be negative.
    This implies that the integrand must vanish a.e. on the support of the measure, i.e.
    \begin{equation*}
        \pi(\tau)v = v
    \end{equation*}
    for all $\tau \in \supp(\mu)$.

    In turn, we have that
    \begin{equation*}
        \operatorname{Re} \langle \hat{\mathcal{L}}(\pi) v, v \rangle =\underbrace{\sum_{j,k=1}^d a_{jk} \langle d\pi(X_j)\,d\pi(X_k) v, v \rangle}_{\leq 0} + \underbrace{\int_{G \setminus \{e\}} \left( \operatorname{Re}\langle \pi(\tau) v, v \rangle - \langle v, v \rangle \right) \mu(d\tau)}_{\leq 0}\leq 0.
    \end{equation*}

    \step{Step 3.}
    We now examine $\hat{\mathcal{L}}(\pi)v = 0$.
    Define $H_0$ to be the closed subgroup generated by the support of $\mu$ and the one-parameter subgroup $\{\exp(tY_l) \mid \forall l,  t \in \mathbb{R}\}$. From step $2$, we see that given $v$ invariant under the closed subgroup $H_0$, the real parts vanish leaving

    \begin{equation*}
        \begin{aligned}
            \hat{\mathcal{L}}(\pi)v & =  \sum_{i,j=1}^d a_{ij}\, d\pi(X_i)\, d\pi(X_j)v + \sum_{i=1}^d b_i\, d\pi(X_i)v                                           \\ & + \int_{G \setminus \{e\}} \left( \pi(\tau)v - v - \sum_{i=1}^d x_i(\tau)\, d\pi(X_i)v \right) \mu(d\tau),\\
                                    & = \sum_{i=1}^d b_i\, d\pi(X_i)v - \int_{G \setminus \{e\}} \left( \sum_{i=1}^d x_i(\tau)\, d\pi(X_i) \right) v\,\mu(d\tau), \\
                                    & = 0.
        \end{aligned}
    \end{equation*}
    Here, we would like to define the Lie algebra element
    \begin{equation}\label{X_b}X_{\tilde{b}} := \sum_{i=1}^d b_i X_i - \int_{G \setminus \{e\}} \left( \sum_{i=1}^d x_i(\tau) X_i \right) \mu(d\tau),\end{equation}
    and conclude that any vector $v$ in the kernel of the symbol must be $d\pi(X_{\tilde{b}})$ invariant. However the integral on the right hand side of \eqref{X_b} is not necessarily convergent as the Lévy measure need only integrate $\lvert x\rvert^2$ in a neighborhood of the identity. Here we have a term that looks like $\lvert x\rvert$. Away from the identity, $x(\tau)$ is bounded so the integral is finite there.

    We thus instead define
    \begin{equation}\label{eqn:drift}X_{\tilde{b}} := \sum_{i=1}^d b_i X_i - \int_{G \setminus U} \left( \sum_{i=1}^d x_i(\tau) X_i \right) \mu(d\tau),\end{equation}
    where $U$ is a neighborhood of the identity specified below, and show that in any direction $v \in \ker(\hat{\mathcal{L}}(\pi))$ , the integrand vanishes.\\

    By the closed subgroup theorem \cite[Thm 3.42.]{MR3331229}, there exists a neighborhood $W$ of $0$ in $\Lie(G)$ such that $H_0 \cap \exp(W) = \exp(\Lie(H_0) \cap W)$ and $\exp$ is a diffeomorphism on $W$. Choosing some $\tau \in H_0 \cap \exp(W)$, we can then write $\log(\tau) \in \Lie(H_0)$, where $\exp(t \log(\tau))$ is an element of $H_0$ for all $t \in \mathbb{R}$.

    The coordinate functions $x_i$ are canonical, so for $\tau \in \supp(\mu) \subset H_0$ near the identity we have $\sum_{i=1}^d x_i(\tau)X_i = \log(\tau) \in \Lie(H_0)$. So,
    \begin{equation*}
        \int_{U} \left( \sum_{i=1}^d x_i(\tau) X_i \right) \mu(d\tau)= \int_{U} \log(\tau)  \mu(d\tau),\end{equation*}

    and $\exp\big(t\sum_{i=1}^d x_i(\tau)X_i\big) \in H_0$ for all $t \in \mathbb{R}$. In turn,
    \begin{equation*}
        d\pi\Big(\sum_{i=1}^d x_i(\tau)X_i\Big)v = 0 \iff \pi\Big(\exp\Big(t\sum_{i=1}^d x_i(\tau)X_i\Big)\Big)v = v.
    \end{equation*}

    Choosing $U \subset \exp(W)$ we then have
    \begin{equation*}
        \int_{U} \Big\lVert \Big( \sum_{i=1}^d x_i(\tau)\, d\pi(X_i) \Big) v \Big\rVert\, \mu(d\tau) = 0.
    \end{equation*}

    We then can conclude that
    \begin{equation*}
        \begin{aligned}
            d\pi(X_{\tilde{b}})v & = \sum_{i=1}^d b_i d\pi(X_i)v - \int_{G \setminus U} \left( \sum_{i=1}^d x_i(\tau) d\pi(X_i) \right)v \mu(d\tau),     \\
                                 & = \sum_{i=1}^d b_i d\pi(X_i)v - \int_{G \setminus \{e\}} \left( \sum_{i=1}^d x_i(\tau) d\pi(X_i) \right)v \mu(d\tau),
        \end{aligned}
    \end{equation*}
    as any contribution of $\tau \in U$ will vanish in the direction of $v.$

    \step{Step 4.}
    The definition~\eqref{eqn:drift} of $X_{\tilde{b}}$ is clearly dependent on the choice of $U \subset \exp(W)$. Without loss of generality assume a different choice of neighborhood  $U'$ is a subset of $U$, then
    \begin{align*}
         & \int_{G \setminus U'} \left( \sum_{i=1}^d x_i(\tau) X_i \right) \mu(d\tau) - \int_{G \setminus U} \left( \sum_{i=1}^d x_i(\tau) X_i \right) \mu(d\tau) \\=
         & \int_{U\setminus U'} \left( \sum_{i=1}^d x_i(\tau) X_i \right) \mu(d\tau).
    \end{align*}
    See then that as the integral is finite and $\sum_{i=1}^d x_i(\tau) X_i \in \Lie(H_0)$ then the integral defines an element  $Z \in \Lie(H_0)$. The subgroup $H$ is closed so $\Lie(H)$ is a subalgebra of $\Lie(G)$. As $Z, X^{U}_{\tilde{b}} \in \Lie(H)$ then $X^{U'}_{\tilde{b}} = X^{U}_{\tilde{b}} - Z \in \Lie(H).$ By definition $H$ then contains $\{\exp(tX^{U'}_{\tilde{b}}) | t \in \mathbb{R}\}$ so the choice of $U$ does not change $H$.

    \step{Step 5.} Finally we can show $H$ periodicity of $u$.
    Let $R_h$ denote right translation. Then
    \begin{align*}
        \langle R_hu, \varphi \rangle = \langle u, R_{h^{-1}}\varphi \rangle.
    \end{align*}
    Taking the Fourier transform
    \begin{equation*}
        \langle \hat{u},  \widehat{\iota \circ R_{h^{-1}}\varphi}  \rangle.
    \end{equation*}

    First note
    \begin{align*}
            \widehat{\iota \circ R_{h^{-1}}\varphi} & = \int_G  (\iota \circ R_{h^{-1}}\varphi)g \pi(g^{-1}) dg, \\
                                                    & =\int_G  (R_{h^{-1}}\varphi)(g^{-1}) \pi(g^{-1}) dg,       \\
                                                    & =\int_G  \varphi(g^{-1}h^{-1}) \pi(g^{-1}) dg.
        \end{align*}
    Take $g \to h^{-1}g^{-1}$
    \begin{align*}
            \int_G  \varphi(g) \pi(gh) dg
             & =\int_G  \varphi(g) \pi(g) dg \pi(h),                     \\
             & = \int_G  (\iota \circ \varphi)(g^{-1}) \pi(g) dg \pi(h), \\
             & = \int_G  (\iota \circ \varphi)(g) \pi(g^{-1}) dg \pi(h), \\
             & = \widehat{(\iota \circ \varphi)}\pi(h).
        \end{align*}
    Thus
    \begin{align*}
            \langle \hat{u},  \widehat{\iota \circ R_{h^{-1}}\varphi} \rangle & = \sum_{\pi \in \hat{G}} \dim(\pi) \operatorname{Tr}(\hat{u}(\pi) \widehat{\iota \circ R_{h^{-1}}\varphi}(\pi)),  \\
                                                                              & =  \sum_{\pi \in \hat{G}} \dim(\pi) \operatorname{Tr}(\hat{u}(\pi) \widehat{\iota \circ R_{h^{-1}}\varphi}(\pi)), \\
                                                                              & =   \sum_{\pi \in \hat{G}} \dim(\pi) \operatorname{Tr}(\hat{u}(\pi) \widehat{(\iota \circ \varphi)}(\pi)\pi(h)),  \\
                                                                              & =   \sum_{\pi \in \hat{G}} \dim(\pi) \operatorname{Tr}(\pi(h)\hat{u}(\pi)\widehat{(\iota \circ \varphi)}(\pi) ),  \\
                                                                              & =  \sum_{\pi \in \hat{G}}\dim(\pi) \operatorname{Tr}(\hat{u}(\pi) \widehat{(\iota \circ \varphi)}(\pi) ),         \\
                                                                              & = \langle \hat{u},  \widehat{\iota \circ \varphi}  \rangle,
        \end{align*}

    where the second last equality uses $\pi(h)\,\hat u(\pi) = \hat{u}(\pi)$ and inverting the Fourier transform, we determine that $u$ is an $H$-periodic distribution.

    The converse follows as a direct consequence of
    \begin{align*}
        \langle \hat{\mathcal{L}}\hat{u}, \Phi \rangle = \sum_{\pi \in \hat{G}} \dim(\pi) \operatorname{Tr}(\hat{\mathcal{L}}(\pi)\hat{u}(\pi)\Phi(\pi)).
    \end{align*}

    As $\pi(h)\hat{u}(\pi) =\hat{u}(\pi)$, then
    \begin{equation*}
        \pi(h)\hat{u}(\pi)v =\hat{u}(\pi)v.
    \end{equation*}
    This implies that the image of $\hat{u}$ is fixed under the action of $h \in H$. If $v\in V_\pi$ is fixed by $H$, then $d\pi(Y_l)v=0$ for every $l$ and $d\pi(X_{\tilde{b}})v=0$, while $\pi(\tau)v=v$ for every $\tau\in\supp\mu$. Substituting into $\hat{\mathcal{L}}(\pi)$, the diffusion term equals $\sum_l d\pi(Y_l)^2v=0$, the jump term $\int_{G\setminus\{e\}}(\pi(\tau)-I)v\,\mu(d\tau)=0$, and the remaining drift and compensator combine to $d\pi(X_{\tilde{b}})v=0$, the contribution from $U$ vanishing as in Step 3. Hence the $H$-fixed vectors lie in the kernel of $\hat{\mathcal{L}}(\pi)$, and we conclude that
    \begin{equation*}
        \hat{\mathcal{L}}(\pi)\hat{u}(\pi)=0.
    \end{equation*}
    Thus, \begin{equation*}\langle \hat{\mathcal{L}}\hat{u}, \varphi \rangle = 0,\end{equation*} and upon Fourier inversion $\langle \mathcal{L}u, \varphi \rangle = 0$. When $H = G$ we have the Liouville property.

\end{proof}

\subsection{Some observations for the nilpotent case}
\label{sec:polynomials}

It is natural to look for a Fourier-analytic proof of Theorem~\ref{thm:main} where $G$ is a nilpotent Lie group such that  $[G,[G,G]]$ is relatively compact. The right substitution for the Peter--Weyl theorem in this case is the Kirillov orbit method~\cite{MR2069175}. The corresponding
group Fourier transform turns $\mathcal{L}$ into a field of operators $\pi(\mathcal{L})$
indexed by the coadjoint orbits, and the scalar condition becomes the requirement that
$\pi(\mathcal{L})$ be injective on every non-trivial representation (the Rockland
condition). For solutions in $L^2$, or with enough decay, this already gives the
Liouville property through hypoellipticity. Bounded/distributional solutions need much more consideration for the following two reasons
\begin{itemize}
    \item The Fourier transform of a bounded function is an operator-valued distribution
          on $\widehat G$, so one needs the full distributional calculus of the orbit method,
    \item The inverse of $\pi(\mathcal{L})$ is not uniformly bounded as the orbits collapse to the
          trivial representation, which is exactly the near-abelian regime where the
          step-dependent behavior discussed above lives.
\end{itemize}
Making this precise would take a development well beyond the scope of the present paper.

Nevertheless, we highlight one subtlety specific to the nilpotent setting. On the Heisenberg group $\mathbb H^1$ write
$\mathcal{L}_0=X^2+Y^2$ for the sub-Laplacian and fix $s\in(0,1)$. There are two standard definitions of the
fractional sub-Laplacian $(-\mathcal{L}_0)^s$. The \emph{Bochner} definition subordinates the heat semigroup
$P_\tau=e^{\tau\mathcal{L}_0}$,
\begin{equation*}
    (-\mathcal{L}_0)^s f=c_s\int_0^\infty\big(f-P_\tau f\big)\,\frac{d\tau}{\tau^{1+s}},\qquad c_s=\frac{s}{\Gamma(1-s)}>0,
\end{equation*}
and makes sense only when this integral converges in the sense of Bochner. The \emph{Fourier-analytic} definition uses the spectral
resolution $-\mathcal{L}_0=\int_0^\infty\lambda\,dE_\lambda$ and sets $(-\mathcal{L}_0)^s=\int_0^\infty\lambda^s\,dE_\lambda$.
Both of these definitions agree on $L^2(\mathbb H^1)$. However, on more general distributions the two definitions can disagree. We demonstrate this on polynomials.
If $f$ is a polynomial with $\mathcal{L}_0 f\neq0$, then
$f-P_\tau f=-\sum_{k\ge1}\tfrac{\tau^k}{k!}\mathcal{L}_0^k f$ is a nonzero polynomial in $\tau$, so the integrand
in the Bochner formula behaves like $\tau^{k-1-s}$ and the integral diverges. Hence, $(-\mathcal{L}_0)^s f$ is
undefined in the Bochner sense. The Fourier-analytic definition, by contrast, remains well defined on such
$f$, by analytic continuation in $s$~\cite{CorniFerrari}. In short, polynomials lie in the domain of the
Fourier-analytic fractional sub-Laplacian but not of the Bochner one. They therefore fall outside the
bounded setting of Theorem~\ref{thm:main}, and can be used only in the Fourier-analytic, distributional
formulation. For the authors, this raises a few questions with regard to the correct notion of solution for the distributional case. This is a future research problem that the authors would like to revisit and understand in more detail.

\section*{Acknowledgements}
The authors would like to thank Anthony Dooley for helpful discussions.

\section*{Declaration on the use of generative AI}
During the preparation of this manuscript the authors used the generative AI assistant Claude (Anthropic), which brought the Choquet--Deny theorem of B\'enard and Frisch~\cite{BenardFrisch} to the authors' attention as the key ingredient for the nilpotent case, and which assisted with editing and exposition. The authors reviewed and verified all content and take full responsibility for it. No generative AI tool is credited as an author.

\end{document}